\documentclass[oribibl]{llncs}
\usepackage{amscd,amssymb,amsmath,latexsym,graphicx,epsfig,color,url}
\usepackage[active]{srcltx}
\usepackage[all]{xy}
\usepackage{hyperref}
\hypersetup{breaklinks=true}

\newcommand{\W}{\mathbb{W}}
\newcommand{\K}{\mathbb{K}}
\newcommand{\V}{\mathbb{V}}

\newcommand{\N}{\mathbb{N}}
\renewcommand{\P}{\mathbb{P}}

\newcommand{\cB}{{\mathcal B}}
\newcommand{\cC}{{\mathcal C}}

\newcommand{\cK}{{\mathcal K}}
\newcommand{\cL}{{\mathcal L}}

\newcommand{\cP}{{\mathcal P}}
\newcommand{\cQ}{{\mathcal Q}}

\newcommand{\be}{{\boldsymbol{e}}}

\newtheorem{ack}[equation]{Acknowledgments}

\begin{document}
\title{Moving curve ideals of rational plane parametrizations}

\author{Carlos D'Andrea}
\institute{ Facultat de Matem\`atiques, Universitat de Barcelona.
Gran Via 585, 08007 Barcelona, Spain
\email{cdandrea@ub.edu \, http://atlas.mat.ub.es/personals/dandrea}
\thanks{Partially supported by the Research Project MTM2010--20279 from the
Ministerio de Ciencia e Innovaci\'on, Spain}}


\maketitle

\begin{abstract}
In the nineties, several methods for dealing in a more efficient way with the implicitization of rational parametrizations were explored in the Computer Aided Geometric Design Community. The analysis of the validity of these techniques has been a fruitful ground for  Commutative Algebraists and Algebraic Geometers, and several results have been obtained so far. Yet, a lot of research is still being done currently around this topic.  In this note we present these methods, show their mathematical formulation, and  survey current results and open questions.
\end{abstract}

\section{Rational Plane Curves}\label{s1}
Rational curves are fundamental tools in Computer Aided Geometric Design. They are used to trace the boundary of any kind of  shape via transforming a parameter (a number) via some simple algebraic operations into a point of the cartesian plane or three-dimensional space. Precision and esthetics in Computer Graphics demands more and more sophisticated calculations, and hence any kind of simplification of the very large list of tasks that need to be performed between the input and the output is highly appreciated in this world. In this survey, we will focus on a simplification of a method for implicitization rational curves and surfaces defined parametrically. This method was developed in the 90's by Thomas Sederberg and his collaborators (see \cite{STD94,SC95,SGD97}), and turned out to become a very rich and fruitful area of interaction among mathematicians, engineers and computer scientist. As we will see at the end of the survey, it is still a very active of research these days.
\begin{figure}[htbp]
\centerline{\includegraphics[angle=0,scale=0.35]{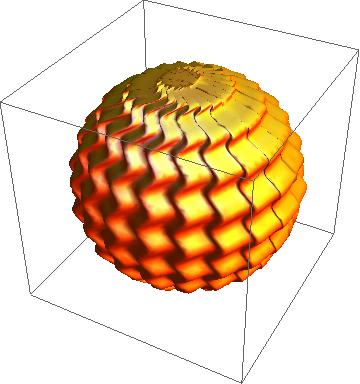}}
 \vspace{3mm}\caption{The shape of an ``orange'' plotted with {\tt Mathematica} 8.0 (\cite{math}).}\label{fig:0}
\end{figure}
\par\smallskip To ease the presentation of the topic, we will work here only with plane curves and point to the reader to the references for the general cases (spatial curves and rational hypersurfaces).

\par\smallskip Let $\K$ be a field, which we will suppose to be algebraically closed so our geometric statements are easier to describe. Here, when we mean ``geometric'' we refer to Algebraic Geometry and not Euclidean Geometry which is the natural domain in Computer Design. Our assumption on $\K$ may look somehow strange in this context, but we do this for the ease of our presentation. We assume the reader also to be familiar with  projective lines and planes over $\K$, which will be denoted with $\P^1$ and $\P^2$ respectively. A {\em rational plane parametrization} is a map
\begin{equation}\label{param}
\begin{array}{cccc}
\phi: & \P^1 & \longrightarrow & \P^2\\
&(t_0:t_1)&\longmapsto&\big(u_0(t_0,t_1):u_1(t_0,t_1):u_2(t_0,t_1)\big),
\end{array}
\end{equation}
where $u_0(t_0,t_1),\,u_1(t_0,t_1),\,u_2(t_0,t_1)$ are polynomials in $\K[T_0,T_1],$ homogeneous, of the same degree $d\geq1,$ and without common factors. We will call $\cC$ to the image of $\phi,$ and refer to it as {\em the rational plane curve parametrized by $\phi$.}
\par\smallskip This definition may sound a bit artificial for the reader who may be used to look at maps of the form
\begin{equation}\label{parafin}
\begin{array}{ccc}
\K&\dasharrow & \K^2\\
t&\longmapsto& \left(\frac{a(t)}{c(t)},\frac{b(t)}{c(t)}\right),
\end{array}
\end{equation}
with $a(t),\,b(t),\,c(t)\in\K[t]$ without common factors, but it is easy to translate this situation to \eqref{param} by  extending this ``map'' (which actually is not defined on all points of $\K$) to one from $\P^1\to\P^2,$ in a sort of {\em continuous} way. To speak about continuous maps, we need to have a topology on $\K^n$ and/or in $\P^n,$ for $n=1, 2$.  We will endow all these sets with the so-called {\em Zariski topology}, which is the coarsest topology that make polynomial maps as in \eqref{parafin} continuous. 
\par\smallskip
Now it should be clear that there is actually an  advantage in working with projective spaces instead of parametrizations as in \eqref{parafin}: our rational map defined in \eqref{param} is {\em actually} a map, and the translation from $a(t),\,b(t),\,c(t)$ to $u_0(t_0,t_1),\,u_1(t_0,t_1),\,u_2(t_0,t_1)$ is very straightforward. The fact that $\K$ is algebraically closed also comes in our favor, as it can be shown that for parametrizations defined over algebraically closed fields  (see \cite{CLO07} for instance), the curve $\cC$ is actually an {\em algebraic variety} of $\P^2$, i.e. it can be described as the zero set of a finite system of homogeneous polynomial equations in  $\K[X_0,X_1,X_2])$.
\par\smallskip More can be said on the case of $\cC$, the {\em Implicitization's Theorem} in \cite{CLO07} states essentially that 
there exists $F(X_0,X_1,X_2)\in\K[X_0,X_1,X_2],$ homogeneous of degree $D\geq1$, irreducible, such that
$\cC$ is actually the zero set of $F(X_0,X_1,X_2)$ in $\P^2,$ i.e. the system of polynomials equations in this case reduces to one single equation. It can be shown that $F(X_0,X_1,X_2)$ is well-defined up to a nonzero constant in $\K$, and it is called {\em the defining polynomial} of $\cC$. The {\em implicitization problem} consists in computing $F$ having as data the polynomials $u_0,\,u_1,\,u_2$ which are the components of $\phi$ as in \eqref{param}.

\begin{example}\label{circc}
Let $\cC$ be the unit circle with center in the origin $(0,0)$ of $\K^2.$ A well-known parametrization of this curve by using a pencil of lines centered in $(-1,0)$ is given in {\em affine} format  \eqref{parafin} as follows:
\begin{equation}\label{pitagoras}
\begin{array}{ccc}
\K&\dasharrow & \K^2\\
t&\longmapsto &  \left(\frac{1-t^2}{1+t^2},\,\frac{2t}{1+t^2},\right).
\end{array}
\end{equation}
Note that if $\K$ has square roots of $-1,$  these values do not belong to the field of definition of the parametrization above. Moreover, it is straightforward to check that the point $(-1,0)$ is not in the image of \eqref{pitagoras}.
\begin{figure}[htbp]
\centerline{\includegraphics[angle=0,scale=0.35]{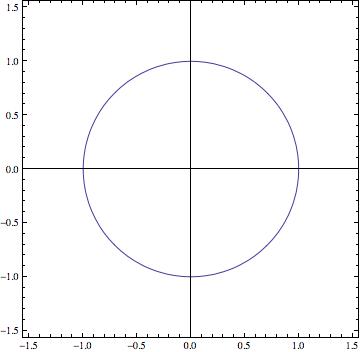}}
 \vspace{-3mm}\caption{The unit circle.}\label{fig:11}
\end{figure}
However, by converting \eqref{pitagoras} into the homogeneous version \eqref{param}, we obtain the parametrization
\begin{equation}\label{circh}
\begin{array}{cccc}
\phi:&\P^1&\longrightarrow & \P^2\\
& (t_0:t_1)& \longmapsto & \big(t_0^2+t_1^2:t_0^2-t_1^2: 2t_0t_1\big),
\end{array}
\end{equation}
which is well defined on all $\P^1.$ Moreover, every point of the circle (in projective coordinates) is in the image of $\phi$, for instance $(1:-1:0)=\phi(0:1),$ which is the point in $\cC$ we were ``missing'' from the parametrization \eqref{pitagoras}. The defining polynomial of $\cC$ in this case is clearly $F(X_0,X_1,X_2)=X_1^2+X_2^2-X_0^2.$
\end{example}
 
In general, the solution to the implicitization problem involves tools from {\em Elimination Theory}, as explained in \cite{CLO07}: from the equation
$$(X_0:X_1:X_2)=\big(u_0(t_0:t_1):u_1(t_0:t_1):u_2(t_0:t_1)\big),$$
one  ``eliminates'' the variables $t_0$ and $t_1$ to get an expression involving only the $X$'s variables.  
\par\smallskip
The elimination process can be done with several tools. The most popular and general is provided by {\em Gr\"obner bases}, as explained in \cite{AL94} (see also \cite{CLO07}).  In the case of a rational parametrization like the one we are handling here, we can consider a more efficient and suitable tool: the {\em Sylvester resultant} of two homogeneous polynomials in $t_0,\,t_1,$ as defined in \cite{AJ06} (see also \cite{CLO05}). We will denote this resultant with $\mbox{\rm Res}_{t_0,t_1}(\cdot,\cdot).$ The following result can be deduced straightforwardly from the section of Elimination and Implicitization in \cite{CLO07}.
\begin{proposition}
There exist  $\alpha,\,\beta\in\N$ such that -up to a nonzero constant-
\begin{equation}\label{resultant}
\mbox{\rm Res}_{t_0,t_1}\big(
X_2u_0(t_0,t_1)-X_0u_2(t_0,t_1),X_2u_1(t_0,t_1)-X_1u_2(t_0,t_1)\big)=X_2^\alpha F(X_0,X_1,X_2)^\beta.
\end{equation}
\end{proposition}
Note that as the polynomial $F(X_0,X_1,X_2)$ is well-defined up to a nonzero constant, all formulae involving it must also hold  this way. For instance,  an explicit computation of \eqref{result} in Example \ref{circc} shows that this resultant is equal to  
\begin{equation}\label{ant}
-4X_2^2\big(X_0^2-X_1^2-X_2^2\big).
\end{equation}
One may think that the number $-4$ which appears above is just a random constant, but indeed it is indicating us something very important: if the characteristic of $\K$ is $2,$ then it is easy to verify that \eqref{pitagoras} does not describe a circle, but the line $X_2=0.$ What is even worse, \eqref{circh} is not the parametrization of a curve, as its image is just the point $(1:1:0).$ 
\par
 To compute the Sylvester Resultant one can use the well-known {\em Sylvester matrix} (see \cite{AJ06,CLO07}), whose nonzero entries contain  coefficients of the two polynomials $X_2u_0(t_0,t_1)-X_0u_2(t_0,t_1)$ and $X_2u_1(t_0,t_1)-X_1u_2(t_0,t_1),$ regarded as polynomials in the variables $t_0$ and $t_1$.  The resultant is then the determinant of that (square) matrix. 

\par\smallskip For instance, in Example \ref{circc}, we have 
 $$\begin{array}{lcl}
X_2\,u_0(t_0,t_1)-X_0\,u_2(t_0,t_1)&=&X_2t_0^2-2X_0t_0t_1+X_2t_1^2\\ X_2\,u_1(t_0,t_1)-X_1\,u_2(t_0,t_1)&=&X_2t_0^2-2X_1t_0t_1-X_2t_1^2,
\end{array}$$
and \eqref{ant} is obtained as the determinant of the Sylvester matrix
\begin{equation}\label{ssyl}
\left(\begin{array}{rrrr}
X_2&-2X_0&X_2&0\\
0&X_2&-2X_0 &X_2\\
X_2&-2X_1&-X_2&0\\
0&X_2&-2X_1&-X_2
\end{array}
\right).
\end{equation}
\par Having $X_2$ as a factor in \eqref{resultant} is explained by the fact that the polynomials whose resultant is being computed in \eqref{result} are not completely symmetric in the $X$'s parameters, and indeed $X_2$ is the only $X$-monomial appearing in both expansions.
\par\smallskip
 The exponent $\beta$ in \eqref{resultant} has a more subtle explanation, it is the {\em tracing index} of the map $\phi$, or the cardinality of its {\em generic fiber}. Geometrically, for all but a finite number of points  $(p_0:p_1:p_2)\in\cC,\,\beta$ is the cardinality of the set $\phi^{-1}(p_0:p_1:p_2).$ Algebraically, it is defined as the degree of the extension $$\left[\K\big(u_0(t_0,t_1)/u_2(t_0,t_1), u_1(t_0,t_1)/u_2(t_0,t_1)\big):\,\K(t_0/t_1)\right].$$
In the applications, one already starts with a map $\phi$ as in \eqref{param} which is {\em generically injective}, i.e. with $\beta=1.$ This assumption is not a big one, due to the fact that generic parametrizations are generically injective, and moreover,  thanks to {\em Lur\"oth's theorem} (see  \cite{vdw66}),  every parametrization  $\phi$  as in \eqref{param} can be factorized as
$\phi=\overline{\phi}\circ\,{\mathcal P},$ with $\overline{\phi}:\P^1\to\P^2$ generically injective, and ${\mathcal P}:\P^1\to\P^1$ being a map defined by a pair of coprime homogeneous polynomial both of them having degree $\beta.$ One can then regard $\overline{\phi}$ as a ``reparametrization'' of $\cC$, and there are very efficient algorithms to deal with this problem, see for instance \cite{SWP08}.
\par\smallskip In closing this section, we should mention the difference between ``algebraic (plane) curves'' and  the  rational curves introduced above. An algebraic plane curve is a subset of $\P^2$ defined by the zero set of a homogeneous polynomial $G(X_0,X_1,X_2)$. In this sense, any rational plane curve is algebraic, as we can find its defining equation via the implicitization described above.
But not all algebraic curve is rational, and moreover, if the curve has degree $3$ or more, a generic algebraic curve will not be rational. Being rational or not is actually a geometric property of the curve,and one should not expect to detect it from the form of the defining polynomial, see \cite{SWP08} for algorithms to decide whether a given polynomial $G(X_0,X_1,X_2)$  defines a rational curve or not. 
\begin{figure}[htbp]
\centerline{\includegraphics[angle=0,scale=0.35]{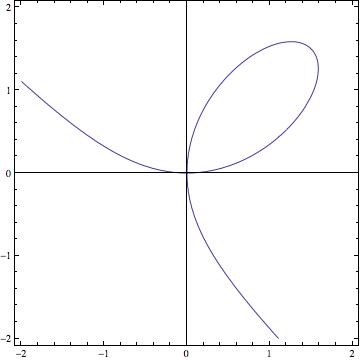}}
 \vspace{-3mm}\caption{The Folium of Descartes.}\label{fig:1}
\end{figure}
For instance, the Folium of Descartes (see Figure \ref{fig:1}) is a rational curve with parametrization
$$(t_0:t_1)\mapsto(t_0^3+t_1^3:3t_0^2t_1:3t_0t_1^2),
$$ and implicit equation given by the polynomial $F(X_0,X_1,X_2)=X_1^3+X_2^3-3X_0X_1X_2.$ On the other hand, Fermat's cubic plotted in Figure \ref{fig:2} is defined by the vanishing of $G(X_0,X_1,X_2)=X_1^3+X_2^3-X_0^3$ but it is not rational.
\begin{figure}[htbp]
\centerline{\includegraphics[angle=0,scale=0.35]{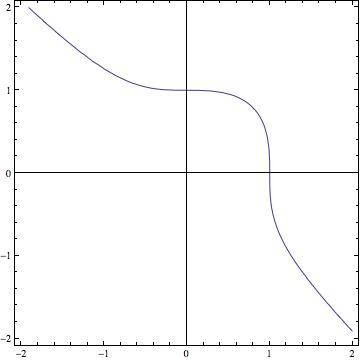}}
 \vspace{-3mm}\caption{Fermat's cubic.}\label{fig:2}
\end{figure}
\par\smallskip
The reason why rational curves play a central role in Visualization and Computer Design  should be easy to get, as they are 
\begin{itemize}
\item easy to ``manipulate''  and be plotted,
\item enough to describe all possible kind of shape by using patches (so-called spline curves).
\end{itemize}

\medskip
\section{Moving lines and $\mu$-bases}
{\em Moving lines} were introduced by Thomas W. Sederberg and his collaborators in the nineties, \cite{STD94,SC95,SGD97,CSC98}. The idea is the following: in each row of the Sylvester matrix appearing in \eqref{ssyl} one can find  the  coefficients as a polynomial in $t_0,\,t_1$ of a form ${\mathcal L}(t_0,t_1,X_0,X_1,X_2)\in\K[t_0,t_1,X_0,X_1,X_2]$ of degree $3$ in the variables $t$'s, and satisfying: 
\begin{equation}\label{ideen}
{\mathcal L}\big(t_0,t_1,u_0(t_0,t_1), u_1(t_0,t_1), u_2(t_0,t_1)\big)=0.
\end{equation}
The first row of \eqref{ssyl} for instance, contains the coefficients of 
$$t_0(X_2\,u_0(t_0,t_1)-X_0\,u_2(t_0,t_1))={\bf X_2} t_0^3{-\bf 2X_0}t_0^2t_1 +{\bf X_2}t_0t_1^2+{\bf0}t_1^3,$$ which clearly vanishes if we set $X_i\mapsto u_i(t_0,t_1).$
Note that all the elements in \eqref{ssyl} are linear in the $X$'s variables.
\par
With this interpretation in mind, we can regard  any such ${\mathcal L}(t_0,t_1,X_0,X_1,X_2)$ as a family of lines in $\P^2$ in such a way that for any $(t_0:t_1)\in\P^1,$ this line passes through the point $\phi(t_0:t_1)\in\cC.$ Motivated by this idea, the following central object in this story has been defined.
\begin{definition}\label{mline}
A {\em moving line} of degree $\delta$ which follows the parametrization $\phi$ is a polynomial
$${\mathcal L}_\delta (t_0,t_1,X_0,X_1,X_2)=v_0(t_0,t_1)X_0+v_1(t_0,t_1)X_1+v_2(t_0,t_1)X_2\in\K[t_0,t_1,X_0,X_1,X_2],$$
with each $v_i$ homogeneous of degree $\delta,\,i=0,1,2$, such that $${\mathcal L}_\delta (t_0,t_1,u_0(t_0,t_1),u_1(t_0,t_1),u_2(t_0,t_1))=0,$$ i.e.
\begin{equation}\label{ecuacion}
v_0(t_0,t_1)u_0(t_0,t_1)+v_1(t_0,t_1)u_1(t_0,t_1)+v_2(t_0,t_1)u_2(t_0,t_1)=0.
\end{equation}
\end{definition}
Note that both $X_2u_0(t_0,t_1)-X_0u_2(t_0,t_1)$ and $X_2u_1(t_0,t_1)-X_1u_2(t_0,t_1)$ are always moving lines following $\phi$. Moreover, note that if we multiply any given moving line by a homogeneous polynomial in $\K[t_0,t_1]$, we obtain another moving line of higher degree. The set of moving lines following a given parametrization  has an algebraic  structure of a {\em module} over the ring $\K[t_0, t_1].$ Indeed, another way of saying that ${\mathcal L}_\delta(t_0,t_1,X_0,X_1,X_2)$ is a moving line which follows $\phi$ is that the vector 
$(v_0(t_0,t_1), v_1(t_0,t_1), v_2(t_0,t_1))$ is a homogeneous element of the {\em syzygy module} of the ideal generated by the sequence $\{ u_0(t_0, t_1), \, u_1(t_0, t_1), \, u_2(t_0, t_1)\}$ -the coordinates of $\phi$- in the ring of polynomials $\K[t_0,t_1].$
\par\smallskip
We will not go further in this direction yet, as the definition of moving lines does not require understanding concepts like syzygies or modules. Note that computing moving lines is very easy from  an equality like \eqref{ecuacion}. Indeed, one first fixes $\delta$ as small as possible, and then sets $v_0(t_0,t_1),\,v_1(t_0,t_1),\,v_2(t_0,t_1)$ as homogeneous polynomials of degree $\delta$ and unknown coefficients, which can be solved via the linear system of equations determined by \eqref{ecuacion}.
\par\smallskip With this very simple but useful object, the {\em method of implitization by moving lines} as stated in \cite{STD94} says essentially the following: look for a set of moving lines of the same degree $\delta,$ with $\delta$ as small as possible, which are ``independent'' in the sense that the matrix of their coefficients (as polynomials in $t_0,\,t_1$) has maximal rank. If you are lucky enough, you will find $\delta+1$ of these forms, and hence the matrix will be square. Compute then the determinant of this matrix, and you will get a non-trivial multiple of the implicit equation. If your are even luckier, your determinant will be equal to $F(X_0,X_1,X_2)^\beta.$

\smallskip
\begin{example}\label{ejj}
Let us go back to the parametrization of the unit circle given in Example \ref{circc}. We check straightforwardly that both
$$\begin{array}{cclcl}
\cL_1(t_0,t_1,X_0,X_1,X_2)&=&-t_1X_0-t_1X_1+t_0X_2&=&X_2\,t_0-(X_0+X_1)\,t_1\\
\cL_2(t_0,t_1,X_0,X_1,X_2)&=&-t_0X_0+t_0X_1+t_1X_2&=&(-X_0+X_1)\,t_0+X_2\,t_1.
\end{array}
$$
satisfy \eqref{ideen}. Hence, they are moving lines of degree $1$ which follow the parametrization of the unit circle. Here, $\delta=1.$ We compute the matrix of their coefficients as polynomials (actually, linear forms) in $t_0,t_1$, and get
\begin{equation}\label{halfmatrix}
\left(\begin{array}{lr}
X_2&-X_0-X_1\\
-X_0+X_1& X_2
\end{array}\right).
\end{equation}
It is easy to check that the determinant of this matrix is equal to $$F(X_0,X_1,X_2)=X_1^2+X_2^2-X_0^2.$$  Note that the size of \eqref{halfmatrix} is actually half of the size of \eqref{ssyl}, and also that the determinant of this matrix gives the implicit equation without any extraneous factor.
\end{example}
\begin{figure}[htpb]
\centerline{  \begin{tabular}{ccc}
\includegraphics[scale=0.35]{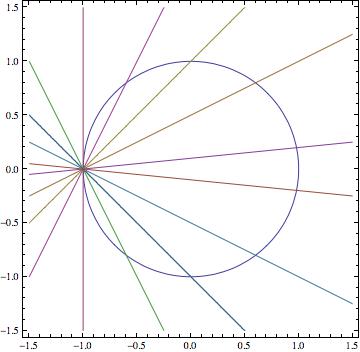}&\hspace*{6mm}& \includegraphics[scale=0.35]{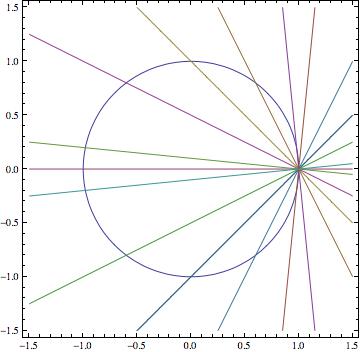}    
  \end{tabular}}
 \vspace{-3mm}\caption{Moving lines $\cL_1$ (left) and $\cL_2$ (right).}\label{fig:6}
\end{figure}
\smallskip
Of course,  in order to convince the reader that this method is actually better than just performing  \eqref{resultant}, we must shed some light on how to compute algorithmically a matrix of moving lines. 
 The following result was somehow discovered by Hilbert more than a hundred years ago, and rediscovered in the CAGD community in the late nineties (see \cite{CSC98}).

\begin{theorem}\label{syz}
For $\phi$ as in \eqref{param}, there exist a unique $\mu\leq\frac{d}{2}$ and two moving lines following $\phi$ which we will denote as  $\cP_\mu(t_0,t_1,X_0,X_1,X_2),\,\cQ_{d-\mu}(t_0,t_1,X_0,X_1,X_2)$ of degrees $\mu$ and $d-\mu$ respectively such that any other moving line following $\phi$ is a polynomial combination of these two, i.e.\ if every $\cL_\delta(t_0,t_1,X_0,X_1,X_2)$  as in the Definition \ref{mline} can be written as
$$\cL_\delta(t_0,t_1,X_0,X_1,X_2)=p(t_0,t_1)\cP_\mu(t_0,t_1,X_0,X_1,X_2)+q(t_0,t_1)\cP_{d-\mu}(t_0,t_1,X_0,X_1,X_2),
$$
with $p(t_0,t_1),\,q(t_0,t_1)\in\K[t_0,t_1]$ homogeneous of degrees $\delta-\mu$ and $\delta-d+\mu$ respectively.
\end{theorem}

This statement is consequence of a stronger one, which essentially says that a parametrization $\phi$ as in \eqref{param}, can be ``factorized'' as follows:
\begin{theorem}[Hilbert-Burch]\label{HB1}
For $\phi$ as in \eqref{param}, there exist a unique $\mu\leq\frac{d}{2}$ and two parametrizations $\varphi_\mu,\,\psi_{d-\mu}:\P^1\to\P^2$ of degrees $\mu$ and $d-\mu$ respectively such that
\begin{equation}\label{wedge}
\phi(t_0:t_1)=\varphi_\mu(t_0:t_1)\times\psi_{d-\mu}(t_0:t_1),
\end{equation}
where $\times$ denotes the usual cross product of vectors.
\end{theorem}
Note that we made an abuse of notation in the statement of \eqref{wedge}, as $\varphi_\mu(t_0:t_1)$ and $\psi_{d-\mu}(t_0:t_1)$ are elements in $\P^2$ and the cross product is not defined in this space. The meaning of $\times$ in \eqref{wedge} should be understood as follows: {\em pick representatives in $\K^3$ of both $\varphi_\mu(t_0:t_1)$ and $\psi_{d-\mu}(t_0:t_1),$ compute the cross product of these two representatives, and then ``projectivize'' the result to $\P^2$ again.}
\par\smallskip The parametrizations $\varphi_\mu$ and $\psi_{d-\mu}$ can be explicited by computing a {\em free resolution} of the ideal 
$\langle u_0(t_0,t_1),\,u_1(t_0,t_1),\,u_2(t_0,t_1)\rangle\subset\K[t_0,t_1],$ and there are algorithms to do that, see for instance \cite{CDNR97}. Note that even though general algorithms for computing free resolutions are based on computations of Gr\"obner bases, which have in general bad complexity time, the advantage here is that we are working with a graded resolution, and also that the resolution of an ideal like the one we deal with here  is of {\em Hilbert-Burch} type  in the sense of \cite{eis95}. This means that  the coordinates of both $\varphi_d$ and $\psi_{d-\mu}$ appear in the columns of the $2\times3$ matrix of the first syzygies in  the resolution. We refer the reader to \cite{CSC98} for more details on the proofs of Theorems \ref{syz} and \ref{HB1}.
\par\smallskip
The connection between the moving lines $\cP_\mu(t_0,t_1,X_0,X_1,X_2),\,\cQ_{d-\mu}(t_0,t_1,X_0,X_1,X_2)$ of Theorem \eqref{syz} and the parametrizations $\varphi_\mu,\,\psi_{d-\mu}$ in \eqref{wedge}  is the obvious one: the coordinates of $\varphi_\mu$ (resp. $\psi_{d-\mu}$) are the coefficients of $\cP_{\mu}(t_0,t_1,X_0,X_1,X_2)$ (resp. $\cQ_{d-\mu}(t_0,t_1,X_0,X_1,X_2)$) as a polynomial in $X_0,\,X_1,\,X_2.$

\begin{definition}
A sequence $\{\cP_\mu(t_0,t_1,X_0,X_1,X_2),\,\cQ_{d-\mu}(t_0,t_1,X_0,X_1,X_2)\}$ as in Theorem \ref{syz}, is called a {\em $\mu$-basis} of $\phi.$ 
\end{definition}

Note that both theorems \ref{syz} and \ref{HB1} only state the uniquenes of the value of $\mu$, and not of $\cP_\mu(t_0,t_1,X_0,X_1,X_2)$ and $\cQ_{d-\mu}(t_0,t_1,X_0,X_1,X_2).$ Indeed, if $\mu=d-\mu$ (which happens generically if $d$ is even), then any two generic linear combinations of the elements of a $\mu$-basis is again another $\mu-$basis.  If $\mu<d-\mu,$ then any polynomial multiple of $\cP_\mu(t_0,t_1,X_0,X_1,X_2)$ of the proper degree can be added to $\cQ_{d-\mu}(t_0,t_1,X_0,X_1,X_2)$ to produce a different $\mu$-basis of the same parametrization.

\smallskip
\begin{example}
For the parametrization of the unit circle given in Example \ref{circc}, one can easily check that
$$\begin{array}{lcl}
\varphi_1(t_0:t_1)&=&(-t_1:-t_1:t_0),\\
\psi_1(t_0:t_1)&=&(-t_0:t_0:t_1)
\end{array}
$$
is a $\mu$-basis of $\phi$ defined in \eqref{circh}, i.e. this parametrization has $\mu=d-\mu=1.$ Indeed, we compute the cross product in \eqref{wedge} as follows: denote with $\be_0,\,\be_1,\,\be_2$ the vectors of the canonical basis of $\K^3$. Then, we get
$$\left|\begin{array}{rrr}
\be_0&\be_1&\be_2\\
-t_1&t_1&t_0\\
-t_0&t_0&t_1
\end{array}
\right|=\big(-t_0^2-t_1^2,\,t_1^2-t_0^2,\,-2t_0t_1\big),
$$
which shows that the $\varphi_1(t_0:t_1)\times\psi_1(t_0:t_1)=\phi(t_0:t_1),$ according to \eqref{wedge}.
\end{example}

The reason the computation of $\mu$-bases is important, is not only because with them we can generate all the moving lines which follow a given parametrization, but also because  they will allow us to produce small matrices of moving lines whose determinant give the implicit equation. Indeed, the following result has been proven in \cite[Theorem 1]{CSC98}.
\begin{theorem}\label{result}
With notation as above, let $\beta$ be the tracing index of $\phi.$ Then, up to a nonzero constant in $\K$, we have
\begin{equation}\label{iimplicit}
\mbox{\rm Res}_{t_0,t_1}\big(\cP_\mu(t_0,t_1,X_0,X_1,X_2),\cQ_{d-\mu}(t_0,t_1,X_0,X_1,X_2)\big)=F(X_0,X_1,X_2)^\beta.
\end{equation}
\end{theorem}

\smallskip
As shown in \cite{SGD97} if you use any kind of matrix formulation for computing the Sylvester resultant, in each row of these matrices, when applied to formulas \eqref{resultant} and \eqref{iimplicit}, you will find the coefficients (as a polynomial in $t_0,\,t_1$) of a moving line following the parametrization. Note that the formula given by Theorem \ref{result} always involves a smaller matrix than the one  in  \eqref{resultant}, as the $t$-degrees of the polynomials $\cP_\mu(t_0,t_1,X_0,X_1,X_2)$ and $\cQ_{d-\mu}(t_0,t_1,X_0,X_1,X_2)$ are roughly half of the degrees of those in \eqref{resultant}.
\par\smallskip
There is, of course, a connection between these two formulas. Indeed,  denote with $\mbox{\rm Syl}_{t_0,t_1}(G,H)$ (resp.\ $\mbox{\rm Bez}_{t_0,t_1}(G,H)$) the Sylvester (resp.\ {\em B\'ezout}) matrix for computing the resultant of two homogeneous polynomials of $G,H\in\K[t_0,t_1].$ For more about definitions and properties of these matrices, see \cite{AJ06}.  In \cite[Proposition 6.1]{BD12}, we prove with Laurent Bus\'e the following:
\begin{theorem}
There exists an invertible matrix $M\in\K^{d\times d}$ such that
$$\begin{array}{l}
X_2\cdot\,\mbox{\rm Sylv}_{t_0,t_1}\big(\cP_\mu(t_0,t_1,X_0,X_1,X_2),\cQ_{d-\mu}(t_0,t_1,X_0,X_1,X_2)\big)\\
=M\cdot\,\mbox{\rm Bez}_{t_0,t_1}\big(
X_2u_0(t_0,t_1)-X_0u_2(t_0,t_1),X_2u_1(t_0,t_1)-X_1u_2(t_0,t_1)\big).
\end{array}
$$
\end{theorem}
From the identity above, one can easily deduce that it is possible to compute the implicit equation (or a power of it) of a rational parametrization with a determinant of a matrix of coefficients of $d$ moving lines, where $d$ is the degree of $\phi$.
Can you do it with less? Unfortunately, the answer is {\it no}, as each row or column of a matrix of moving lines is linear in $X_0, X_1, X_2,$ and the implicit equation has typically degree $d.$ So, the method will work optimally with a matrix of size $d\times d$, and essentially you will be computing the Sylvester matrix of a $\mu$-basis of $\phi$. 

\medskip
\section{Moving conics, moving cubics...}
One can actually take advantage of the resultant formulation given in \eqref{iimplicit} and get a determinantal formula for the implicit equation by using the square matrix $$\mbox{\rm Bez}_{t_0,t_1}\big(\cP_\mu(t_0,t_1,X_0,X_1,X_2),\cQ_{d-\mu}(t_0,t_1,X_0,X_1,X_2)\big),$$ which has smaller size (it will have $d-\mu$ rows and columns) than the Sylvester matrix of these polynomials.  But this will not be a matrix of coefficients of moving lines anymore, as the input coefficients of the B\'ezout matrix will be quadratic in $X_0,\,X_1,\,X_2.$ Yet, due to the way the B\'ezout matrix is being built (see for instance \cite{SGD97}, one can find in the rows of this matrix the coefficients of a polynomial which also vanishes on the parametrization $\phi.$ This motivates the following definition:

\begin{definition}\label{mcurve}
A {\em moving curve} of bidegree $(\nu,\delta)$ which follows the parametrization $\phi$ is a polynomial ${\mathcal L}_{\nu,\delta}(t_0,t_1,X_0,X_1,X_2)
\in\K[t_0,t_1,X_0,X_1,X_2]$ homogeneous in $X_0,\,X_1,\,X_2$ of degree $\nu$ and in $t_0,\,t_1$ of degree $\delta,$ such that 
$${\mathcal L}\big(t_0,t_1,u_0(t_0,t_1),u_1(t_0,t_1),u_2(t_0,t_1)\big)=0.$$ 
\end{definition}
If $\nu=1$ we recover the definition of moving lines given in \eqref{mline}.  For $\nu=2,$ the polynomial $\cL(t_0,t_1,X_0,X_1,X_2)$ is called a {\em moving conic} which follows $\phi$ (\cite{ZCG99}). {\em Moving cubics} will be curves with $\nu=3,$ and so on.
\par A series of experiments made by Sederberg and his collaborators showed something interesting: one can compute the defining polynomial of $\cC$ as a determinant of a matrix of coefficients of moving curves following the parametrization, but the more singular the curve is (i.e.\ the more singular points it has), the smaller the determinant of moving curves gets. For instance, the following result appears in \cite{SC95}:
\begin{theorem}\label{sc95}
The implicit equation of a quartic curve with no base points can be written as a $2\times2$ determinant. If the curve doesn't have a triple point, then each element of the determinant is a quadratic; otherwise one row is linear and one row is cubic.
\end{theorem}
To illustrate this, we consider the following examples.
\begin{example}\label{exx}
Set $u_0(t_0,t_1)=t_0^4-t_1^4,\,u_1(t_0,t_1)=-t_0^2t_1^2, u_2(t_0,t_1)=t_0t_1^3.$ These polynomials defined a parametrization $\phi$ as in \eqref{param} with implicit equation given by the polynomial  $F(X_0,X_1,X_2)=X_2^4- X_1^4 - X_0X_1X_2^2.$ From the shape of this polynomial, it is easy to show that $(1:0:0)\in\P^2$ is a point of multiplicity $3$ of this curve, see Figure \ref{fig1}. In this case, we have $\mu=1,$ and it is also easy to verify that
$$\cL_{1,1}(t_0,t_1,X_0,X_1,X_2)=t_0X_2+t_1X_1$$ is a moving line which follows $\phi$. The reader will now easily check that the following moving curve of bidegree $(3,1)$ also follows $\phi$:
$$\cL_{1,3}(t_0,t_1,X_0,X_1,X_2)=t_0(X_1^3+X_0X_2^2)+t_1\,X_2^3.
$$ 
And the $2\times2$ matrix claimed in Theorem \ref{sc95} for this case is made  with the coefficients of both $\cL_{1,1}(t_0,t_1,X_0,X_1,X_2)$ and $\cL_{1,3}(t_0,t_1,X_0,X_1,X_2)$ as polynomials in $t_0,\,t_1:$
$$\left(\begin{array}{lr}
X_2&X_1\\
X_1^3+X_0X_2^2&X_2^3
\end{array}\right).
$$
\begin{figure}[htbp]
\centerline{\includegraphics[angle=0,scale=0.35]{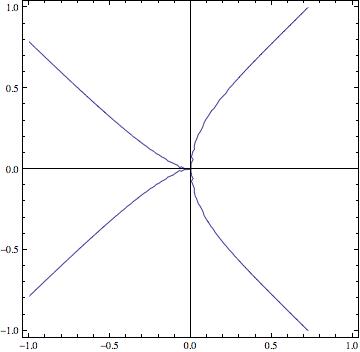}}
 \vspace{-3mm}\caption{The curve of Example \ref{exx}.}\label{fig1}
\end{figure}
\end{example}

\smallskip
\begin{example}\label{xex}
We reproduce here Example 2.7 in \cite{cox08}. Consider
$$u_0(t_0,t_1)=t_0^4,\,u_1(t_0,t_1)=6t_0^2t_1^2-4t_1^4,\, u_2(t_0,t_1)=4t_0^3t_1-4t_0t_1^3.$$
This input defines a quartic curve with three nodes, with implicit equation given by \newline $F(X_0,X_1,X_2)=X_2^4+4X_0X_1^3+2X_0X_1X_2^2-16X_0^2X_1^2-6X_0^2X_2^2+16X_0^3X_1,$ see Figure \ref{fig111}.
\par\smallskip
The following two moving conics of degree $1$ in $t_0,\,t_1$ follow the parametrization:
$$\begin{array}{ccl}
\cL_{1,2}(t_0,t_1,X_0,X_1,X_2)&=&t_0(X_1X_2-X_0X_2)+t_1(-X_2^2-2X_0X_1+4X_0^2)\\
\tilde{\cL}_{1,2}(t_0,t_1,X_0,X_1,X_2)&=& t_0(X_1^2+\frac12X_2^2-2X_0X_1)+t_1(X_0X_2-X_1X_2).
\end{array}
$$ As in the previous example, the $2\times2$ matrix of the coefficients of these moving conics is the matrix claimed in Theorem \ref{sc95}.
\begin{figure}[htbp]
\centerline{\includegraphics[angle=0,scale=0.35]{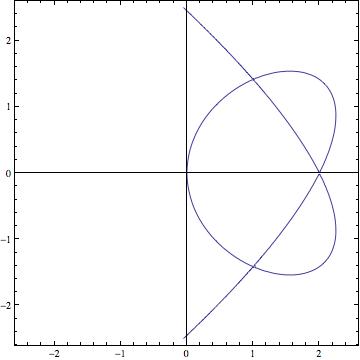}}
 \vspace{-3mm}\caption{The curve of Example \ref{xex}.}\label{fig111}
\end{figure}
\end{example}

\medskip
\section{The moving curve ideal of $\phi$}
Now it is time to introduce some tools from Algebra which will help us understand all the geometric constructions defined above.The set of all moving curves following a given parametrization generates a {\em bi-homogeneous} ideal in $\K[t_0,t_1,X_0,X_1,X_2]$, which we will call the {\em moving curve ideal} of this parametrization.
\par\smallskip
As explained above,  the method of moving curves for implicitization of a rational parametrization looks for small determinants made with coefficients of moving curves which follow the parametrization of low degree in $t_0,\,t_1.$ To do this, one would like to have a description as in Theorem \ref{syz}, of a set of ``minimal'' moving curves from which we can  describe in an easy way all the other elements of the moving curve ideal.
\par \smallskip
Fortunately, Commutative Algebra provides the adequate language and tools for dealing with this problem. As it was shown by David Cox in \cite{cox08}, all we have to do is look for minimal generators of the kernel $\cK$ of the following morphism of rings:
\begin{equation}\label{rees}
\begin{array}{cclr}
\K[t_0,t_1,X_0,X_1,X_2]&\longrightarrow&\K[t_0,t_1,z]\\
t_i&\longmapsto& t_i&\,i=0,1,\\
X_j&\longmapsto& u_j(t_0,t_1)\,z&\,j=0,1,2.
\end{array}
\end{equation}
Here, $z$ is a new variable. The following result appears in \cite[Nice Fact 2.4]{cox08} (see also \cite{BJ03} for the case when $\phi$ is not generically injective):
\begin{theorem}
$\cK$ is the moving curve ideal of $\phi$ .
\end{theorem}
Let us say some words about the map \eqref{rees}. Denote with $I\subset\K[t_0,\,t_1]$ the ideal generated by $u_0(t_0,t_1),\,u_1(t_0,t_1),\,u_2(t_0,t_1)$. The image of \eqref{rees} is actually isomorphic to $\K[t_0,t_1][z\,I],$ which is called the {\em Rees Algebra} of $I$. By the Isomorphism Theorem, we then get that $\K[t_0,t_1,X_0,X_1,X_2]/\cK$ is isomorphic to the Rees Algebra of $I$. This is why the generators of $\cK$ are called the {\em defining equations} of the Rees Algebra of $I$. The Rees Algebra that appears in the moving lines method corresponds to the blow-up of $V(I),$ the variety defined by $I$. Geometrically, it is is just the blow-up of the empty space (the effect of this blow-up is just to introduce torsion...),  but yet the construction should explain somehow why moving curves are sensitive to the presence of complicated singularities. It is somehow strange that the fact that the description of $\cK$ actually gets much simpler if the singularities of $\cC$ are more entangled.
\par\smallskip
Let us show this with an example. It has been shown in \cite{bus09}, by unravelling some duality theory developed by Jouanolou in \cite{jou97}, that for any proper parametrization of a curve of degree $d$  having $\mu=2$ and only cusps as singular points, the kernel $\cK$ has $\frac{(d+1)(d-4)}{2}+5$ minimal generators. On the other hand, in a joint work with Teresa Cortadellas \cite{CD13b} (see also \cite{KPU13}), we have shown that if $\mu=2$ and there is a point of very high multiplicity (it can be proven that if the multiplicity of a point is larger than $3$ when $\mu=2$, then it must be equal to $d-2$), then the number of generators drops to $\lfloor\frac{d+6}2\rfloor,$ i.e. the description of $\cK$ is simpler in this case.  In both cases, these generators can be make explicit, see \cite{bus09, CD13b, KPU13}.
\par\smallskip
Further evidence supporting this claim is what is already known for the case $\mu=1$, which was one of the first one being worked out by several authors: \cite{HSV08,CHW08,bus09,CD10}. It turns out (cf.\ \cite[Corollary 2.2]{CD10}) that $\mu=1$ if and only if the parametrization is proper (i.e. generically injective), and there is a point on $\cC$ which has multiplicity $d-1$, which is the maximal multiplicity a point can have on a curve of degree $d.$ If this is the case, then the parametrization has exactly $d+1$ elements. 
\par\smallskip In both cases ($\mu=1$ and $\mu=2$), explicit elements of a set of minimal generators of $\cK$ can be given in terms of the input parametrization. But in general, very little is known about how many are them and which are their bidegrees. 
Let $n_0(\cK)$ be the $0$-th {\em Betti number} of $\cK$ (i.e.\  the cardinal of any minimal set of generators of $\cK$). We propose the following problem which is the subject of attention of several researchers at the moment.

\begin{problem}\label{1prob}
Describe  {\em all} the possible values of $n_0(\cK)$ and the parameters that this function depends on, for a proper parametrization $\phi$ as in \eqref{param}.
\end{problem}
Recall that ``proper'' here means ``generically injective''. For instance, we just have shown above that, for $\mu=1,\,n_0(\mu)=d+1.$ If $\mu=2,$ the value of $n_0(\cK)$ depends on whether there is a very singular point or not. Is $n_0$ a function of only $d,\,\mu$ and the multiplicity structure of $\cC$?  
\par
A more ambitious problem of course is the following. Let $\cB(\cK)\subset\N^2$ be the (multi)-set of bidegrees of a minimal set of generators of $\cK$.
\begin{problem}\label{2prob}
Describe {\em all} the possible values of $\cB(\cK).$
\end{problem}
For instance, if $\mu=1,$ we have that (see \cite[Theorem 2.9]{CD10})
$$\cB(\cK)=\{(0,d),\,(1,1),\,(1,d-1),\,(2,d-2),\ldots,(d-1,1)\}.
$$
Explicit descriptions of $\cB(\cK)$ have been done also for $\mu=2$ in \cite{bus09,CD13b,KPU13}.  In this case, the value of $\cB(\cK)$ depends on whether the parametrization has singular point of multiplicity $d-2$ or not.
\par For $\mu=3$ the situation gets a bit more complicated as we have found in \cite{CD13b}:
consider the parametrizations $\phi_1$ and $\phi_2$ whose $\mu$-bases are respectively :
$$\begin{array}{l}
\cP_{3,1}(t_0,t_1,X_0,X_1,X_2)=t_0^3X_0+(t_1^3-t_0t_1^2)X_1\\
\cQ_{7,1}(t_0,t_1,X_0,X_1,X_2)=(t_0^6t_1-t_0^2t_1^5)X_0+(t_0^4t_1^3+t_0^2t_1^5)X_1+(t_0^7+t_1^7)X_2,
\\ \\
\cP_{3,2}(t_0,t_1,X_0,X_1,X_2)=(t_0^3-t_0^2t_1)X_0+(t_1^3+t_0t_1^2-t_0^2t_1)X_1 \\ 
\cQ_{7,2}(t_0,t_1,X_0,X_1,X_2)=(t_0^6t_1-t_0^2t_1^5)X_0+(t_0^4t_1^3+t_0^2t_1^5)X_1+(t_0^7+t_1^7)X_2.
\end{array}$$
Each of them parametrizes properly a rational plane curve of degree $10$ having the point $(0:0:1)$ with multiplicity 7. The rest of them are either double or triple points. Set $\cK_1$ and $\cK_2$ for the respective kernels, we have then
$$\begin{array}{rcl}
\cB(\cK_1)&=&\{(3,1),\,(7,1),\,(2,3),\,(2,3),\,(4,2),\,(2,4),\,(1,6),\,(1,6),\,(1,6),\,(0,10)\},\\
\cB(\cK_2)&=&\{(3,1),\,(7,1),\,(2,3),\,(2,3),\,(4,2),\,(2,4),\,{\bf (1,5)},\,(1,6),\,(1,6),\,(0,10)\}.
\end{array}$$
The parameters to find in the description of $n_0(\cK)$ proposed in Problem \ref{1prob} may be more than $\mu$ and the multiplicities of the curve. For instance, in \cite{CD13}, we have shown that if there is a minimal generator of bidegree $(1,2)$ in $\cK,$ then the whole 
set $\cB(\cK)$ is constant, and equal to
$$\left\{\begin{array}{lrc}
\{(0,d),\,(1,2),\,(1,d-2),\,(2,d-4),\ldots,(\frac{d-1}{2},1),\,(\frac{d+1}2,1)\}&\,&\mbox{if}\,d\,\mbox{is odd}\\
\{(0,d),\,(1,2),\,(1,d-2),\,(2,d-4),\ldots,(\frac{d}{2},1),\,(\frac{d}2,1)\}&\,&\mbox{if}\,d\,\mbox{is even.}
\end{array}
\right.
$$
To put the two problems above in a more formal context, we proceed as in \cite[Section 3]{CSC98}: For $d\geq1,$ denote with $\V_d\subset{\K[t_0,t_1]_d}^3$ the set  of triples of homogeneous polynomials $(u_0(t_0,t_1),\,u_1(t_0,t_1),\,u_2(t_0,t_1)\big)$ defining a proper parametrization $\phi$ as in \eqref{param}. Note that one can regard $\V_d$ as an open set in an algebraic variety in the space of parameters. Moreover,  $\V_d$ could be actually be taken as a quotient of ${\K[t_0,t_1]_d}^3$ via the action of $\mbox{SL}(2,\K)$ acting on the monomials $t_0,\,t_1$.
\begin{problem}\label{3prob}
Describe the subsets of $\V_d$  where $\cB(\cK)$ is constant.
\end{problem}
Note that, naturally the $\mu$-basis is contained in $\cK,$ and moreover, we have (see \cite[Proposition 3.6]{BJ03}):
$$\cK= \langle \cP_\mu(t_0,t_1,X_0,X_1,X_2),\,\cQ_{d-\mu}(t_0,t_1,X_0,X_1,X_2)\rangle:\langle t_0,\,t_1\rangle^\infty,
$$
so the role of the $\mu$-basis is crucial to understand $\cK$. Indeed, any minimal set of generators of $\cK$ contains a $\mu$-basis, so the pairs $(1,\mu),\,(1,d-\mu)$ are always elements of $B(\cK).$ The study of the geometry of $\V_d$ according to the stratification done by $\mu$ has been done in  \cite[Section 3]{CSC98} (see also \cite{dan04,iar13}).  Also, in \cite{CKPU13}, a very interesting study of how the $\mu$-basis of a parametrization having generic $\mu$ ($\mu=\lfloor d/2\rfloor$) and very singular points  looks like has been made.  It would be interesting to have similar results for $\cK$. 
\par\smallskip
In this context, one could give a positive answer to the experimental evidence provided by Sederberg and his collaborators about the fact that ``the more singular the curve, the simpler the description of $K$'' as follows. For $\W\subset\V_d$, we denote by $\overline{\W}$ the closure of $\W$ with respect to the Zariski topology.

\begin{conjecture}\label{1conj}
If $\W_1,\,\W_2\subset\V_d$ are such that $n_0|_{\W_i}$ is constant for $i=1,2,$ and $\overline{\W}_1\subset\overline{\W}_2,$ then
$$n_0\big(\W_1\big)\leq n_0\big(\W_2\big).$$
\end{conjecture}
Note that this condition is equivalent to the fact that $n_0(\cK)$ is {\em upper semi-continuous} on $\V_d$ with its Zariski topology. Very related to this conjecture is the following claim, which essentially asserts that in the ``generic'' case, we obtain the largest value of $n_0(\cK):$
\begin{conjecture}\label{2conj}
Let $\W_d$ be open set of $\V_d$ parametrizing all the curves with   $\mu=\lfloor d/2\rfloor,$ and having all its singular points being ordinary double points. Then, $n_0(\cK)$ is constant on $\W_d,$ and attains its maximal value on $\V_d$ in this component.
\end{conjecture}
Note that a ``refinement'' of Conjecture \ref{1conj} with $\cB(\cK_1)\subset\cB(\cK_2)$ will not hold in general, as in the examples computed for $\mu=2$ in \cite{bus09,CD13b,KPU13} show. Indeed, we have in this case that the Zariski closure of those parametrizations with a point of multiplicity $d-2$ is contained in the case where all the points are only cusps, but the bidegrees of the minimal generators of $\cK$ in the case of parametrizations with points of multiplicity $d-2$ appear at lower values than the more general case (only cusps).

\medskip
\section{Why  Rational Plane Curves only?}
All along this text we were working with the parametrization of a rational plane curve, but most of the concepts, methods and properties worked out here can be extended in two different directions. The obvious one is to consider ``surface'' parametrizations, that is maps of the form
\begin{equation}\label{paramS}
\begin{array}{cccc}
\phi_S:&\P^2&\dasharrow&\P^3\\
&(t_0:t_1:t_2)&\longmapsto&\big(u_0(t_0,t_1,t_2):u_1(t_0,t_1,t_2):u_2(t_0,t_1,t_2):u_3(t_0,t_1,t_2)\big)
\end{array}
\end{equation}
where $u_i(t_0,t_1,t_2)\in\K[t_0,t_1,t_2],\,i=0,1,2,3,$ are homogeneous of the same degree, and without common factors. Obviously, one can do this in higher dimensions also, but we will restrict the presentation just to this case. The reason we have now a dashed arrow  in \eqref{paramS} is because even with the conditions imposed upon the $u_i$'s, the map may not be defined on all points of $\P^2.$ For instance, if $$u_0(t_0,t_1,t_2)=t_1t_2,\,u_1(t_0,t_1,t_2)=t_0t_2,\,u_2(t_0,t_1,t_2)=t_0t_1,\, u_3(t_0,t_1,t_2)=t_0t_1+t_1t_2,$$ $\phi_S$ will not be defined on the set $\{(1:0:0),\,(0:1:0),\,(0:0:1)\}.$
\par\smallskip
In this context, there are methods to deal with the implicitization analogues to those presented here for plane curves. For instance, one can use a {\em multivariate resultant} or a {\em sparse resultant} (as defined in \cite{CLO05}) to compute the implicit equation of the Zariski closure of the image of $\phi_S$. Other tools from Elimination Theory such as determinants of complexes can be also used to produce matrices whose determinant (or quotient or $\gcd$ of some determinants) can also be applied to compute the implicit equation, see for instance \cite{BJ03,BCJ09}.
\par\smallskip
The method of moving lines and curves presented before gets translated into a {\em method of moving planes and surfaces} which follows $\phi_S$, and its description and validity is much more complicated, as the both the Algebra and the Geometry involved have more subtleties,  see \cite{SC95,CGZ00,cox01,BCD03,KD06}.
Even though it has been shown in \cite{CCL05} that there exists an equivalent of a $\mu$-basis in this context, its computation of is not as easy as in the planar case. Part of the reason is that the syzygy module of general $u_i(t_0,t_1,t_2),\,i=0,1,2,3$ is not free anymore (i.e. it does not have sense the meaning of a ``basis'' as we defined in the case of curves), but if one set $t_0=1$ and regards these polynomials as affine bivariate forms, a nicer situation appears but without control on the degrees of the elements of the $\mu$-basis, see \cite[Proposition 2.1]{CCL05} for more on this. Some explicit descriptions have been done for either low degree parametrizations, and also for surfaces having some additional geometric features (see \cite{CSD07,WC12,SG12,SWG12}), but the general case remains yet to be explored.
\par \smallskip
A generalization of a map like \eqref{rees} to this situation is straightforward, and one can then consider the defining ideal of the Rees Algebra associated to $\phi_S$. Very little seems to be known about the minimal generators of $\cK$ in this situation. In \cite{CD10} we studied the case of {\em monoid} surfaces, which are rational parametrizations with a point of the highest possible multiplicity. This situation can be regarded as a possible generalization of the case $\mu=1$ for plane curves, and has been actually generalized to {\em de Jonqui\`eres} parametrizations in \cite{HS12}.
\par \smallskip
We also dealt in  \cite{CD10} (see also \cite{HW10}) with the case where there are two linearly independent moving planes of degree $1$ following the parametrization plus some geometric conditions, this may be regarded of a generalization of the ``$\mu=1$'' situation for plane curves. But the general description of the defining ideal of the Rees Algebra for the surface situation is still an open an fertile area for research.
\par\smallskip
The other direction where we can go after consider rational plane parametrizations is to look at spatial curves, that is maps
$$\begin{array}{cccc}
\phi_C:&\P^1&\longrightarrow&\P^3\\
&(t_0:t_1)&\longmapsto&(u_0(t_0,t_1):u_1(t_0,t_1):u_2(t_0,t_1):u_3(t_0,t_1)),
\end{array}
$$
where  $u_i\in\K[t_0,t_1],$ homogeneous of the same degree $d\geq1$ in $\K[t_0,t_1]$ without any common factor. In this case, the image of $\phi_C$ is a curve in $\P^3$, and one has to replace ``an'' implicit equation with ``the'' implicit equations, as there will be more than one in the same way that the implicit equations of the line joining $(1:0:0:1)$ and $(0:0:0:1)$ in $\P^3$ is given by the vanishing of the equations $X_1=X_2=0.$
\par\smallskip
As explained in \cite{CSC98}, both Theorems \ref{syz} and \ref{HB1} carry on to this situation, so there is more ground to play and theoretical tools to help with the computations. In \cite{CKPU13}, for instance, the singularities of the spatial curve are studied as a function of the shape of the $\mu$-basis. Further computations have been done in \cite{KPU09} to explore the generalization of the case $\mu=1$ and produce generators for $\cK$ in this case. These generators, however, are far from being minimal.  More explorations have been done in \cite{JG09,HWJG10,JWG10}, for some specific values of the degrees of the generators of the $\mu$-basis.
\par\smallskip 
It should be also mentioned that in the recently paper \cite{iar13}, an attempt of the stratification proposed in Problem \ref{2prob} for this kind of curves is done, but only with respect to the the value of $\mu$ and no further parameters.
\par\smallskip
As the reader can see, there are lots of recent work in this area, and many many challenges yet to solve. We hope that in the near future we can get more and deeper insight in all these matters, and also to be able to apply these results in the Computer Aided and Visualization community.

\medskip
\begin{ack}
I am grateful to Laurent Bus\'e, Eduardo Casas-Alvero and Teresa Cortadellas Benitez for their careful reading of a preliminary version of this manuscript, and very helpful comments. Also, I thank the anonymous referee for her further comments and suggestions for improvements, and to Marta Narv\'aez Clauss for her help with the computations of some examples. All the plots in this text have been done with {\tt Mathematica 8.0} (\cite{math}).
\end{ack}

\end{document}